\documentclass[a4paper,12pt]{amsart}
\usepackage{amssymb}
\usepackage{ifthen}
\usepackage{graphicx}
\usepackage{float}
\usepackage{caption}
\usepackage{subcaption}
\usepackage{cite}
\usepackage{amsfonts}
\usepackage{amscd}
\usepackage{amsxtra}
\usepackage{color}

\newtheorem{thm}{Theorem}
\newtheorem{cor}{Corollary}
\newtheorem{lem}{Lemma}
\newtheorem{rem}{Remark}

\newtheorem{claim}{Claim}

\newtheorem{conj}{Conjecture}
\newtheorem{prob}{Problem}

\theoremstyle{definition}
\newtheorem{defn}{Definition}[section]
\newtheorem{example}{Example}

\newenvironment{pf}[1][]{%
 \vskip 1mm
 \noindent
 \ifthenelse{\equal{#1}{}}%
  {{\slshape Proof. }}%
  {{\slshape #1.} }%
 }%
{\qed\bigskip}

\newcounter{alphabet}

\newenvironment{Thm}[1][]{\refstepcounter{alphabet}%
\bigskip%
\noindent%
{\bf Theorem \Alph{alphabet}}%
\ifthenelse{\equal{#1}{}}{}{ (#1)}%
{\bf .} \itshape}{\vskip 8pt}
\newenvironment{Lem}[1][]{\refstepcounter{alphabet}%
\bigskip%
\noindent%
{\bf Lemma \Alph{alphabet}}%
{\bf .} \itshape}{\vskip 8pt}

\newcommand{\T}{{\mathbb T}}

\newcommand{\IR}{{\mathbb R}}

\newcommand{\IC}{{\mathbb C}}
\newcommand{\ID}{{\mathbb D}}

\newcommand{\IZ}{{\mathbb Z}}

\newcommand{\D}{{\mathbb D}}

\def\be{\begin{equation}}
\def\ee{\end{equation}}

\newcommand{\bee}{\begin{enumerate}}
\newcommand{\eee}{\end{enumerate}}

\newcommand{\blem}{\begin{lem}}
\newcommand{\elem}{\end{lem}}
\newcommand{\bthm}{\begin{thm}}
\newcommand{\ethm}{\end{thm}}
\newcommand{\bcor}{\begin{cor}}
\newcommand{\ecor}{\end{cor}}
\newcommand{\beg}{\begin{example}}
\newcommand{\eeg}{\end{example}}
\newcommand{\begs}{\begin{examples}}
\newcommand{\eegs}{\end{examples}}
\newcommand{\bdefe}{\begin{defn}}
\newcommand{\edefe}{\end{defn}}
\newcommand{\bprob}{\begin{prob}}
\newcommand{\eprob}{\end{prob}}
\newcommand{\bques}{\begin{ques}}
\newcommand{\eques}{\end{ques}}
\newcommand{\bei}{\begin{itemize}}
\newcommand{\eei}{\end{itemize}}

\newcommand{\bde}{\begin{deter}}
\newcommand{\ede}{\end{deter}}
\newcommand{\bca}{\begin{case}}
\newcommand{\eca}{\end{case}}
\newcommand{\bcl}{\begin{claim}}
\newcommand{\ecl}{\end{claim}}
\newcommand{\bcon}{\begin{conj}}
\newcommand{\econ}{\end{conj}}
\newcommand{\bcons}{\begin{conjs}}
\newcommand{\econs}{\end{conjs}}
\newcommand{\bprop}{\begin{propo}}
\newcommand{\eprop}{\end{propo}}
\newcommand{\br}{\begin{rem}}
\newcommand{\er}{\end{rem}}
\newcommand{\brs}{\begin{rems}}
\newcommand{\ers}{\end{rems}}
\newcommand{\bo}{\begin{obser}}
\newcommand{\eo}{\end{obser}}
\newcommand{\bos}{\begin{obsers}}
\newcommand{\eos}{\end{obsers}}
\newcommand{\bpf}{\begin{pf}}
\newcommand{\epf}{\end{pf}}
\newcommand{\ba}{\begin{array}}
\newcommand{\ea}{\end{array}}
\newcommand{\beq}{\begin{eqnarray}}
\newcommand{\beqq}{\begin{eqnarray*}}
\newcommand{\eeq}{\end{eqnarray}}
\newcommand{\eeqq}{\end{eqnarray*}}

\newcounter{minutes}
\divide\time by 60
\newcounter{hours}
\multiply\time by 60 \addtocounter{minutes}{-\time}

\begin{document}
\title[Relations of the class $\mathcal{U}(\lambda)$ to other families of functions]
{Relations of the class $\mathcal{U}(\lambda)$ to other families of functions}


\thanks{
File:~\jobname .tex,
          printed: \number\day-\number\month-\number\year,
          \thehours.\ifnum\theminutes<10{0}\fi\theminutes}

\author{Liulan Li,   Saminathan  Ponnusamy and Karl-Joachim Wirths
}
\address{L. Li, College of Mathematics and Statistics
 (Hunan Provincial Key Laboratory of Intelligent Information Processing and Application),\\
 Hengyang Normal University,
 Hengyang,  Hunan 421002, People's Republic of China.}
\email{lanlimail2012@sina.cn}

\address{S. Ponnusamy,  Department of Mathematics,
Indian Institute of Technology Madras, Chennai-600 036, India.}
\email{samy@iitm.ac.in}

\address{K.-J. Wirths, Institute of Algebra and Analysis,
Technical University Braunschweig, 38106 Braunschweig, Germany.}
\email{kjwirths@tu-bs.de}





\begin{abstract}
In this article, we consider the family of  functions $f$ analytic in the unit disk
$|z|<1$ with the normalization $f(0)=0=f'(0)-1$ and satisfying the
condition $\big |\big (z/f(z)\big )^{2}f'(z)-1\big |<\lambda $ for
some $0<\lambda \leq 1$. We denote this class by  $\mathcal{U}(\lambda)$ and we are interested in the relations between $\mathcal{U}(\lambda)$ and other families of functions holomorphic or harmonic in the unit disk. Our first example in this direction is the family of functions convex in one direction. Then we are concerned with the subordinates to the function $1/((1-z)(1-\lambda z))$. We prove that not all functions $f(z)/z$ $(f \in \mathcal{U}(\lambda))$ belong to this family. This disproves an assertion from \cite{OPW}. Further, we disprove a related coefficient conjecture for $\mathcal{U}(\lambda)$. We consider the intersection of the class of the above subordinates and $\mathcal{U}(\lambda)$ concerning the boundary behaviour of its functions. At last, with the help of functions from $\mathcal{U}(\lambda)$, we construct functions harmonic and close-to-convex in the unit disk.

\end{abstract}

\maketitle
\pagestyle{myheadings}
\markboth{L. Li, S. Ponnusamy and K.-J. Wirths}{Relations of the class $\mathcal{U}(\lambda)$ to other families of functions  }

\section{Introduction}\label{LSW18-sec1}

 Let $\ID=\{z\in \IC: |z|<1\}$ be the open unit disk and $\mathbb{T}=\{z\in \IC: |z|=1\}$, the unit circle. Let ${\mathcal A}$ be the family of all functions $f$ analytic in $\ID $ with the Taylor series expansion $f(z)=z+\sum_{k=2}^{\infty}a_kz^k$. Let ${\mathcal S}$ denote the subset of ${\mathcal A}$ consisting of functions that
are univalent in $\ID$. See \cite{Duren:univ} for the general
theory of univalent functions.  For $0<\lambda \leq 1$, consider the class
$$ \mathcal{U}(\lambda)=\left \{f\in {\mathcal A}:\, |U_{f}(z)|< \lambda \  \mbox{for}\ z\in \ID \right
\},
$$
where
$$U_{f}(z)=\left(\frac{z}{f(z)}\right)^2f'(z)-1=\frac{z}{f(z)}-z\left(\frac{z}{f(z)}\right)'-1.$$
 A well-known fact about this class is that every $f\in
\mathcal{U}:=\mathcal{U}(1)$ is  univalent in $\ID$ and hence,
$\mathcal{U}(\lambda)\subset \mathcal{U} \subset \mathcal{S}$ for
$\lambda \in (0,1]$  (cf. \cite{Aks58,AksAvh70,OzNu72}). Several
properties of this family are established in \cite{OPW}. For
example, the class ${\mathcal U}(\lambda )$ is preserved under a
number of elementary operations such as rotation, conjugation,
dilation and omitted-value transformations \cite[Lemma~1]{OPW}.
These properties do hold for the family $\mathcal S$. On the other
hand, it was also pointed out in \cite{OPW} that, although $\mathcal
S$ is preserved under the square-root transformation, the family
${\mathcal U}$ (and hence ${\mathcal U}(\lambda )$) does not.
Because $f'(z)\left(z/f(z)\right)^2$ $(f\in {\mathcal U})$ is
bounded, we see that $\left(z/f(z)\right)^{2}f'(z)\neq 0$ in $\ID$
and thus, each $f\in {\mathcal U}$ is  non-vanishing in
$\ID\backslash\{0\}$. It is well-known (cf
\cite{FR-2006,OP-01,PV2005,PW}) that ${\mathcal U}$ neither is
included in ${\mathcal S}^{\star}$ nor includes ${\mathcal
S}^{\star}$. Here ${\mathcal S}^{\star}$ denotes the class of
starlike functions, namely, functions $f\in{\mathcal S}$ such that
$f(\ID)$ is starlike with respect to the origin. Typical members of
$\mathcal{U}$ are in the set ${\mathcal S}_{\IZ}$, where
$${\mathcal S}_{\IZ}=\left \{z, ~~\frac{z}{(1\pm z)^2},~~\frac{z}{1\pm z},~~\frac{z}{1\pm z^2}.
~~ \frac{z}{1\pm z+z^2}\right \}
$$
It is worth pointing out that ${\mathcal
S}_{\IZ}\subset\mathcal{U}\cap \mathcal{S}^*$ and the set ${\mathcal
S}_{\IZ}$  is precisely the set of functions in $\mathcal{S}$ having
integral coefficients in the power series expansion (cf \cite{Fr1}).
Furthermore, each of these functions plays an important role in
function theory since they together with rotations are extremal for
some well-known subfamilies of ${\mathcal S}$.

Let $B$ denote the set of the functions $\omega$ which are analytic
in $\D$ and satisfy $|\omega(z)|\leq 1$ for all $z\in\D$. Also, we
consider the subfamily $B_0=\{\omega\in B:\, \omega(0)=0 \}$. It is
a simple exercise to see that each $f\in{\mathcal U}(\lambda)$ has
the characterization (cf. \cite{obpo-2007a}),
 \be\label{OPW7-eq1}
\frac{z}{f(z)}=1-a_{2}z+\lambda z\int_{0}^{z}\omega(t)\,dt, \ee for
some $\omega\in B$, where $a_2 =f''(0)/2$. By simplification, we may
rewrite \eqref{OPW7-eq1} as
\be\label{LS1}
\frac{z}{f(z)}=1-a_{2}z+\lambda \omega_1(z), \ee where $\omega_1\in
B_0$ and $\omega_1'(0)=0$.

In Section \ref{LSW18-sec2}, we find some functions from ${\mathcal U}(\lambda)$ that are convex in one direction and in addition a member of
${\mathcal U}$ that is convex in no direction.


In \cite[Theorems 3 and 4]{OPW} it was stated that $f\in{\mathcal U}(\lambda )$ satisfies the subordination relations
\be\label{OPW7-eq4c} \frac{f(z)}{z}\prec
\frac{1}{1+(1+\lambda)z+\lambda z^2}=\frac{1}{(1+z)(1+\lambda z)}, ~ z\in \ID,
\ee
and
$$ \frac{z}{f(z)}+ a_2z \prec 1+2\lambda z + \lambda z^2, ~ z\in \ID.
$$
Here $\prec$ denotes the usual subordination \cite{Duren:univ}.
Section \ref{LSW18-sec3} is devoted to these two subordination
relations and to the study of the boundary behaviour of the Schwarz
functions that define subordination. One of the results of this
discussion is the fact that \eqref{OPW7-eq4c} is not valid. The
formula \eqref{OPW7-eq4c} was indeed the motivation for the
following coefficient conjecture of the functions $f(z)=z+
\sum_{n=2}^{\infty}a_n(f)z^k\in {\mathcal U}(\lambda).$

\bcon\label{LiPoWi-conj1}
If $f\in {\mathcal U}(\lambda)$,  then
$$|a_n(f)| \leq \sum_{k=0}^{n-1}\lambda^k \ \mbox{for}\ n\geq 2.
$$
\econ

The conjecture for $n=2$ has been shown in \cite{VY2013} and one of
its alternate proofs was given in \cite{OPW}.
%
One of our aims of this paper is to present a counterexample to this conjecture for the case $n=3$.

In Section \ref{LSW18-sec4} we use functions from the family
$${\mathcal U}_2(\lambda ):= {\mathcal U}(\lambda )\cap\{f\in \mathcal{A}: f''(0)=0\}
$$
and construct functions that are harmonic and close-to-convex (univalent) in $\ID$.

\section{Convex in some direction}\label{LSW18-sec2}

We recall that a univalent mapping $f$ analytic in $\ID$ is called
convex if $f(\ID)$ is a convex domain. A function $f$ analytic in
$\ID$ is said to be close-to-convex if
the complement $D=\IC\backslash f(\ID)$ of $f(\ID)$ can be filled with non-intersecting half-lines
emanating on the boundary $\partial D$ and lying completely in $D$. Every close-to-convex function
is known to be univalent in $\ID$. However, close-to-convexity of harmonic mappings requires univalence in
the definition. This will be considered in Section \ref{LSW18-sec4}.

A domain $D$ in $\IC$ is said to be convex in the direction $\gamma$
$(0\leq \gamma < \pi)$ if every line parallel to the line through
$0$ and $e^{i\gamma}$ has a connected intersection with $D$.  A
univalent function $f$ in $\mathbb{D}$ is said to be convex in the
direction $\gamma$ if $f(\mathbb{D})$ is convex in the direction
$\gamma$.

Obviously, every function that is convex in the direction $\gamma$
is necessarily close-to-convex, but the converse is not true.
Clearly, a convex function is convex in every direction. The class
of functions convex in one direction has been introduced by
Robertson \cite{Rober36} and  investigated by many others (see, for
example, \cite{Lecko02,RZ}). In particular, convex in the direction
$\gamma=0$ is said to be convex along horizontal/real direction.

In order to prove Theorem \ref{convex in a direction}, we need the
following results.

\begin{Thm}\label{convex in some direction theorem}
{\rm (Royster and Ziegler \cite{RZ})} Let $\varphi$ be a
non-constant analytic function in $\ID$. The function $\varphi$ maps
$\ID$  univalently  onto a domain convex in the direction $\gamma $
if and only if there are numbers $\mu$ and $\nu$, $0\leq \mu <2\pi$
and $0\leq \nu\leq \pi$, such that
$$
{\rm Re\,} \{P_{\mu,\,\nu,\,\gamma,\, \varphi}(z)  \}\geq 0,~\mbox{ $z\in \ID$},
$$
where
$$ P_{\mu,\,\nu,\,\gamma,\, \varphi}(z)= e^{i(\mu-\gamma
)}(1-2ze^{-i\mu}\cos \nu+z^2e^{-2i\mu})\varphi'(z).$$
\end{Thm}

\begin{Thm}\label{special form f}{\rm (\cite{OPW, VY2013})}
If $f(z)=z+\sum_{k=2}^{\infty}a_kz^k\in {\mathcal U}(\lambda)$, then
$|a_2|\leq 1+\lambda$. Moreover, if $|a_2|=1+\lambda $, then $f$
must be of the form \beq\label{LiPoWi-Eqn B.1}
f_{\theta}(z)=\frac{z}{1-(1+\lambda)e^{i\theta}z+\lambda
e^{2i\theta}z^2},\,\,\theta\in [0,2\pi). \eeq
\end{Thm}

\br\label{LSW18-rem1}
By a standard argument, a consequence of this theorem is that the functions $f_{\theta}$ are extreme points of the class ${\mathcal U}(\lambda)$. It is an open problem, whether there exist other extreme points of ${\mathcal U}(\lambda)$. See also Remark \ref{LSW18-rem2} below.
\er

\begin{thm}\label{convex in a direction}
 The function $f_{\theta}$ defined in Theorem B is convex in the direction $2\pi-\theta$.
\end{thm}

\bpf  Let $f_{\theta}(z)$ be defined by \eqref{LiPoWi-Eqn B.1}. Then
\be\label{LS18-eq1} f_{\theta}'(z)=\frac{1-\lambda
e^{2i\theta}z^2}{\left(1-\lambda
e^{i\theta}z\right)^2\left(1-e^{i\theta}z\right)^2}. \ee We next
consider $P_{\mu,\, \nu,\, \gamma,\, f_{\theta}}(z)$ for $\mu =
2\pi-\theta$, $\nu=0$ and $\gamma= 2\pi-\theta$ to get
$$ P_{ 2\pi-\theta,\, 0,\, 2\pi-\theta,\, f_{\theta}}(z) = (1-e^{i\theta}z)^2f_{\theta}'(z) = \frac{(1-\lambda e^{2i\theta} z^2)(1-\lambda e^{-i\theta}\overline{z})^2}{|1-\lambda e^{i\theta}z|^4}.$$
This gives by \eqref{LS18-eq1} that
\[
  P_{ 2\pi-\theta,\, 0,\, 2\pi-\theta,\, f_{\theta}}(z)  =\frac{1-\lambda e^{2i\theta}z^2}{\left(1-\lambda e^{i\theta}z\right)^2}
=\frac{N(e^{i\theta}z)}{D(e^{i\theta}z)},
\]
where $ N(z) =\left(1-\lambda z^2\right)\left(1-\lambda
\overline{z}\right)^2$ and $D(z)=|1-\lambda z|^4>0$ in $\ID$.
Clearly, $f_{\theta}$ is convex in the direction $2\pi -\theta$ by
Theorem~A 
if ${\rm Re\,} \{N(z)\}\geq 0$ for $|z|=1$, since 
$$P_{ 2\pi-\theta,\, 0,\,
2\pi-\theta,\, f_{\theta}}(0)=1>0.$$ Thus, it suffices to show that
${\rm Re\,} \{N(z)\}>0$ for $|z|=1$. It is a simple to see that
$${\rm Re\,} \{N(e^{i\theta})\} =1-\lambda^3-2\lambda (1-\lambda)\cos\theta -\lambda (1-\lambda)\cos 2\theta \geq  1-\lambda^3-3\lambda (1-\lambda)=  (1-\lambda)^3
$$
and hence,
$$ {\rm Re\,}\left \{ P_{ 2\pi-\theta,\, 0,\, 2\pi-\theta,\, f_{\theta}}(z)\right \}= {\rm Re\,}\left\{(1-e^{i\theta}z)^2 f_{\theta}'(z)\right\}>0 ~\mbox{ for $z\in\ID$},
$$
showing that $f_{\theta}$ is convex in the direction $2\pi -\theta$ by Theorem~A. 
\epf

Theorem \ref{convex in a direction} motivates us to consider whether
all the functions in ${\mathcal U}(\lambda)$ are convex in some
direction. Next we show that the answer is negative in general.

\begin{figure}[h]
\begin{center}
\includegraphics[height=6.0cm, width=5.5cm, scale=1]{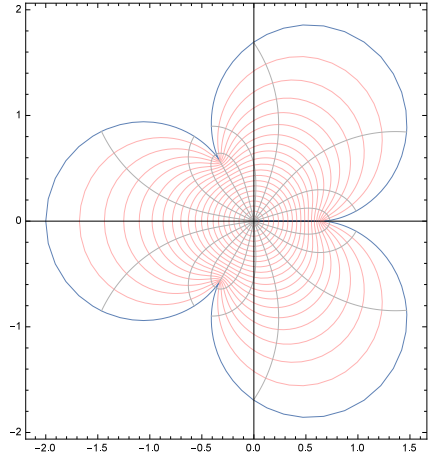}
\end{center}
\caption{The image of unit disk under the function $g(z)= \frac{z}{1+\frac{1}{2}z^3}$\label{LSW-fig1}}
\end{figure}

\begin{example}\label{not convex in any direction}
If $g(z)= \frac{z}{1+\frac{1}{2}z^3}$, then $g$ is not convex in any direction.
 \end{example}
\bpf
First we see that
$$\left(\frac{z}{g(z)}\right)^2g'(z)-1 =-z^2 ~\mbox{ and }~ g(e^{\frac{2i\pi}{3}}z)= e^{\frac{2\pi i}{3}}g(z),
$$
and thus, it follows that $g \in \mathcal{U}$ and $g$ is 3-fold symmetric.
Therefore,
$$D=:g\left(\bigg\{z=re^{i\theta}:\ 0\leq r<1,
0\leq\theta<\frac{2\pi}{3}\bigg\}\right),\;\ g(\ID)=D\bigcup
e^{\frac{2i\pi}{3}}D\bigcup e^{\frac{4i\pi}{3}}D,$$
where $e^{i\theta}D=\{e^{i\theta}z:\,z\in D\}.$ Thus, we only
need to prove that $g$ is not convex in the direction $\gamma$ for
$\gamma\in[0, \frac{\pi}{3})$. There are three inward pointing
(zero-angle) cusps, which are at $z=1,e^{i\frac{2}{3}\pi},
e^{i\frac{4}{3}\pi}$. Therefore, $g$ is not convex in the direction
$\gamma$ for $\gamma\in[0, \frac{\pi}{3})$. See Figure \ref{LSW-fig1}.
 \epf

Naturally one may ask the following:
\bprob\label{LSW18-prob1}
Which functions in ${\mathcal U}(\lambda)$ are convex in some direction?
\eprob

The characterization \eqref{LS1} shows that each $f\in {\mathcal U}(\lambda)$ is completely determined by
its second coefficient $a_2$, $\lambda$ and the Schwarz function
$\omega_1$. One of the simple forms of $f$ is when $a_2=0$ and
$\omega_1(z)=az^n$ ($n\geq2$ and $|a|\leq 1$). So we begin to answer Problem  \ref{LSW18-prob1} in this case.

\begin{thm}\label{convex in some direction 2}
Suppose that $f\in {\mathcal U}_2(\lambda)$ is in the form
\eqref{LS1}, where  $\lambda \in (0,1].$ If $\omega_1(z)=az^2$ and
$|a|\leq1$, then $f$ is convex in the direction $\gamma = 2\pi -
\frac{\arg a}{2}$.
\end{thm}

\bpf There is nothing to prove if $a=0$. So, we assume that $a\neq
0$. To prove the assertions, we take in Theorem A $\mu = \gamma,$
and $\nu = \pi/2.$ We have from \eqref{LS1} with $a_2=0$ that
$$f'(z)=\frac{1-\lambda az^2}{(1+\lambda az^2)^2}=\frac{1-\lambda p s^2}{(1+\lambda p s^2)^2},
$$
where $s=z e^{i\frac{\arg a}{2}}$ and $p=|a|$. By Theorem~A, 
we need to show that for all
$|s|<1$, \be\label{LS2} {\rm Re}\left\{f'(z)(1+s^2)\right\}={\rm
Re}\left\{\frac{(1-\lambda p s^2)(1+s^2)}{(1+\lambda
ps^2)^2}\right\}\geq 0. \ee

If $s=0$, then $f'(z)(1+s^2)$ assumes the value $1$ and thus, it suffices to prove \eqref{LS2} for all $|s|=1$.
Now

\vspace{6pt}
${\rm Re}\left\{f'(z)(1+s^2)\right\}\big|_{|s|=1}$

\beqq
& =& {\rm Re}\left\{\frac{(1-\lambda^2 p^2-\lambda p s^2+\lambda p
\overline{s}^2)(1+\lambda p+s^2+\lambda p\overline{s}^2)}{|1+\lambda ps^2|^4}\right\}\\
&=&(1-\lambda p)\frac{1+4\lambda p+\lambda^2 p^2+(1+\lambda p)^2{\rm Re\,}(s^2)-2\lambda p\left({\rm Re\,}(s^2)\right)^2}{|1+\lambda
ps^2|^4}\geq0
\eeqq
for all $|s|=1$.
 \epf

\begin{figure}[h]
\begin{center}
\includegraphics[height=6.0cm, width=5.5cm, scale=1]{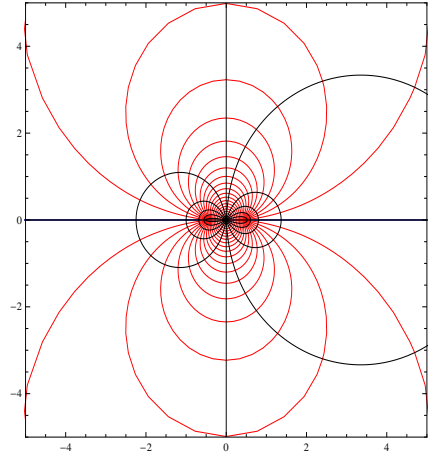}
\hspace{1cm}
\end{center}
\caption{The range $f({\mathbb{D}})$ for $\lambda=1$, $\omega_1(z)= z^2$ \label{LSW-fig2a}
}

%
\end{figure}

We end this section with a picture to demonstrate the image of the
function $f(z)$ involved in Theorem \ref{convex in some direction
2} with certain special values of $\lambda$ and $a$. See Figure \ref{LSW-fig2a}.

\section{Relations between Schwarz functions and  $f\in{\mathcal U}(\lambda)$ }\label{LSW18-sec3}

To prove our next result, we need the following result.

\begin{Thm}\label{LiPoWi*_Thm C}
{\rm (\cite[Proposition 4.13]{P-1992})}
Let $\varphi$ be analytic in $\ID$ with an angular limit $\varphi(\zeta)$ at $\zeta \in \mathbb{T}.$ If $\varphi(\ID) \subset \ID,$ $\varphi(\zeta) \in \mathbb{T},$ then the angular derivative $\varphi'(\zeta)$ exists and
$$0<\zeta\frac{\varphi'(\zeta)}{\varphi(\zeta)}=\sup_{z \in \ID} \frac{1-|z|^2}{|\zeta-z|^2}\frac{|\varphi(\zeta)-\varphi(z)|^2}{1-|\varphi(z)|^2}\leq+\infty.
$$
\end{Thm}

\begin{thm}\label{extreme case}
Let $\frac{z}{f(z)}=1-a_{2}z+\lambda \omega_1(z)\in {\mathcal
U}(\lambda)$, where $\omega_1\in B_0$ and $\omega_1'(0)=0$. If there
exists a point $\zeta$ of the unit circle $\T$ such that
$|\omega_1(\zeta)| =1$, then $\omega_1$ is of the form $\omega_1(z)=
e^{i\theta} z^2$.
\end{thm}

\bpf Let $\omega_1(z)=z^2\psi(z)$. Then $\psi$ is analytic in $\ID$,
$\psi(\ID)\subset\overline{\ID}$, and $|\psi(\zeta)|=1$. If
$\psi(\ID)\subset\ID$, then Theorem~C 
implies that the angular derivative $\psi'(\zeta)$ exists and satisfies that
\beq\label{LiPoWirth*_eqn 3.1}
0<\frac{\zeta\psi'(\zeta)}{\psi(\zeta)}\leq+\infty.
\eeq
 We now
consider
$$U_{f}(z)=-\lambda z^2\left(\psi(z)+z\psi'(z)\right).$$
As $z$ approaches $\zeta$ along the segment connecting $0$ and
$\zeta$, $$|U_{f}(z)|\rightarrow
\lambda\left(1+\frac{\zeta\psi'(\zeta)}{\psi(\zeta)}\right)>\lambda,$$
which is a contradiction. Thus $\psi(z)\equiv e^{i\theta}$ for some
$\theta\in\IR$. The proof is completed.
 \epf

\beg\label{not subordinate}
 Let $\lambda\in(0,1)$, $k\geq 2$ and
$$\frac{z}{f(z)}=(1-z)\left(1-\frac{\lambda z}{k}\sum_{\nu=0}^{k-1}z^{\nu}\right)=1-z\left(1+\frac{\lambda}{k}\right)+\frac{\lambda
z^{k+1}}{k}.$$ Then $f(z)\in {\mathcal U}(\lambda)$ but $z/f(z)$
is not subordinate to $(1-z)(1-\lambda z).$
 \eeg

\bpf Since $(1-z)\left(1-\frac{\lambda
z}{k}\sum_{\nu=0}^{k-1}z^{\nu}\right)$ does not vanish in the unit
disk $\ID$, and
$$
  \frac{z}{f(z)}-z\left(\frac{z}{f(z)}\right)'-1=-\lambda z^{k+1},$$
the function $f$ belongs to ${\mathcal U}(\lambda)$.

Suppose on the contrary that $z/f(z)\prec (1-z)(1-\lambda z).$ Then
there exists a function $\psi\in B$ such that
$$1-(1+\lambda)z\psi(z)+\lambda
z^2\psi^2(z)=1-\left(1+\frac{\lambda}{k}\right)z+\frac{\lambda
z^{k+1}}{k}.$$ Since $k\geq 2$, Theorem~B 
implies that $|\psi(z)|<1$ for $z\in\ID$.
Setting $z=1$ on both sides of the last relation gives
$$(1-\psi(1))(1-\lambda \psi(1))=0.
$$
Moreover, because $\lambda\in(0,1)$, we obtain that $\psi(1)=1$ and
$\psi'(1)=0,$ which contradicts Theorem~C, 
namely \eqref{LiPoWirth*_eqn 3.1}. Therefore, $z/f(z)$ is not subordinate to $(1-z)(1-\lambda z).$ 
\epf

Based on the erroneous subordination statement that if $f(z)\in {\mathcal U}(\lambda)$ then $z/f(z)$ is subordinate to $(1-z)(1-\lambda z)$,
Conjecture \ref{LiPoWi-conj1} was proposed in \cite{OPW}. The idea of Example \ref{not subordinate} was the investigation of functions in ${\mathcal U}(\lambda)$
that have a pole at the point $z=1$. Now such functions will be used to give a counterexample to Conjecture \ref{LiPoWi-conj1}.

\beg\label{a counterexample}
For any $\lambda\in (0,\delta)$, where $\delta=\frac{3\,-\,4\ln  2}{4\ln 2 -\,2} \approx 0.29435$,
there exists $f\in {\mathcal U}(\lambda)$ such that
$a_3(f)\,>\,1+\lambda +\lambda^2.$
\eeg
\bpf
  For $a\in (0,1)$, we consider the functions
  \[ f_a(z)\,=\,\frac{z}{(1-z)\left(1\,-\,\frac{\lambda z}{1-z} \int_z^1\frac{t+a}{1+at}\,dt\right)}\,=
  \,\frac{z}{1\,-\,z\left(1\,+\,\lambda\int_z^1\frac{t+a}{1+at}\,dt\right)}.\]
  Obviously, the functions $f_a$ are holomorphic in $\ID$, and since
  \[\frac{z}{f_a(z)}\,-\,z\left(\frac{z}{f_a(z)}\right)'\,-\,1\,=\,-\lambda z^2 \frac{z+a}{1+az},\]
  they belong to the class ${\mathcal U}(\lambda)$.  Setting
  $$
  v(a):=\int_0^1\frac{t+a}{1+at}\,dt\,=\,\frac{1}{a}\,-\,\frac{1-a^2}{a^2}\ln(1+a)\,=\,\frac{1}{2}\,+\,\sum_{k=1}^{\infty}\frac{2(-1)^{k-1}a^k}{k(k+2)}.
  $$
The Taylor expansion of $f_a$ can be derived from
\beqq
f_a(z)\,&=&\,z\sum_{k=0}^{\infty}z^k\left(1\,+\,\lambda v(a)-\,\lambda \int_0^z\frac{t+a}{1+at}\,dt\right)^k\\
&=& z\,+\,z^2\left(1\,+\,\lambda v(a)\right)\,+\,z^3\left(\left(1\,+\,\lambda v(a)\right)^2\,-\,\lambda a\right)\,+\,o(z^3).
 \eeqq
In particular, the coefficient of $z^3$ of $f_{a}$ gives
  \[a_3(f_a)\,=\,1\,+\,\lambda(2v(a)-a)\,+\,\lambda^2v^2(a),\]
  and we prove that there exist real numbers $a\in (0,1),\, \lambda\in (0,1)$ such that
  \begin{equation}\label{fa}
    1\,+\,\lambda(2v(a)-a)\,+\,\lambda^2v^2(a)\,>\,1+\lambda+\lambda^2.
  \end{equation}
  Firstly, we will be concerned with the derivatives of $w(a):=2v(a)-a$,
  i.e.,
\beqq
w'(a)\,=\,2v'(a)-1\, &=&\frac{2}{a}\,-\,\frac{4}{a^2}\,+\,\frac{4}{a^3}\ln(1+a)-1,\\
&=& 4\sum_{n=3}^{\infty}\frac{(-1)^{n-1} }{n} \, a^{n-3} \, -1,
\eeqq
 so that $w'(0)=\frac{1}{3} ,\ \ w'(1)=4\ln 2\,-\,3\,<\,0 .$ Clearly
$$w''(a)=4\sum_{n=4}^{\infty}\frac{(-1)^{n-1} (n-3) }{n}\, a^{n-4}
$$
so that $w''(0)=-1.$ Moreover, it can be easily seen that
\beqq
w''(a)\,=\,\frac{2}{a^4(1+a)}\left(6a+3a^2-a^3-6(1+a)\ln(1+a)\right):=\frac{2u(a)}{a^4(1+a)}.
  \eeqq

As $u(0)=0,$ $u'(a)=6a-3a^2-6 \ln(1+a)$, $u'(0)=0$, and
\[u''(a)\,=\,\frac{-6a^2}{1+a}\,<\,0,\,a\in(0,1],
\]
the function $u'$ is monotonically decreasing on $[0,1]$ and negative on $(0,1]$, and the function $u$ is negative on $(0,1]$. Hence, $w''(a)\,<\,0$ for $a\in [0,1]$. This implies that the function $w$ is concave on $[0,1]$. Since $w(0)=w(1)=1,$ we get the inequality $w(a)\,>\,1$ for $a\in (0,1).$

  On the other hand,
$$v(a)=\dfrac{w(a)+a}{2},~ 2v'(a)=w'(a)+1 ~\mbox{ and }~ 2v''(a)=w''(a)<0 ~\mbox{ for $a\in [0,1]$},
$$
which implies that   $v'(a)$ is monotonically decreasing for $a\in [0,1]$. So
$$v'(a) >v'(1)=2\ln 2-1>0 ~\mbox{ for $a\in [0,1)$}.
$$
Therefore, $v(a)$ is monotonically increasing for $a\in [0,1]$ and $\frac{1}{2}=v(0)\leq v(a)\leq v(1)=1$.

  The above estimates for $w(a)$ and $v(a)$ yield that
 for any $a\in (0,1)$ and
  \[   0\,<\,\lambda\,<\,\frac{w(a)-1}{1-v^2(a)}\]
 inequality (\ref{fa}) holds. We observe that $\frac{w(a)-1}{1-v^2(a)}$ is monotonically increasing for $a\in [0,1]$.
To get the more precise about the above assertion on the values of $\lambda$ for which (\ref{fa}) is valid for $a\in (0,1)$, we calculate by L'H\^opital's rule
$$  \lim_{a\to 1}\frac{w(a)-1}{1-v^2(a)}\,= \lim_{a\to 1}\frac{w'(a)}{-2v(a)v'(a)}\,=-\frac{w'(1)}{v(1)(w'(1)+1)}=\frac{3\,-\,4\ln 2}{4\ln 2 -\,2}\,.
 $$
    This proves the assertion of our example.
 \epf

\br\label{LSW18-rem2}
Example \ref{a counterexample} shows that for any $\lambda \in (0, \delta),\ \delta=0.2943\cdots$, there are other extreme points in ${\mathcal U}(\lambda)$ besides the functions of the form
$$f_{\theta}(z)=\frac{z}{1-(1+\lambda)e^{i\theta}z+\lambda e^{2i\theta}z^2},\,\,\theta\in [0,2\pi).
$$\er

In the rest of this section, we consider $f\in {\mathcal
U}(\lambda)$ of the form \be\label{PW-eq5} \frac{f(z)}{z}=
\frac{1}{1-(1+\lambda)\phi(z)+\lambda \phi^2(z)}, ~ z\in \ID, \ee
where $\phi\in B_0$. That is, $z/f(z)$ is subordinate to
$(1-z)(1-\lambda z).$

The authors in \cite{PW} discussed the $\phi$ such that $f$ in the
form of \eqref{PW-eq5} belongs to ${\mathcal U}(\lambda)$, which is
represented as follows.

\begin{Thm}\label{impossible condition1}
Let $f\in {\mathcal U}(\lambda)$ be given by \eqref{PW-eq5} with a
function $\phi\in B_0$ analytic on the closed unit disk
$\overline{\ID}$. If there exists a point $\zeta$ of the unit circle
$\T$ such that $\phi(\zeta) = -1$, then $\phi$ is of the form
$\phi(z)= e^{i\theta} z$.
\end{Thm}

Julia's lemma plays a key role in the proof of Theorem~D. 
We will show that Theorem~D 
still holds even if $\phi$ is not analytic
on the whole closed unit disk $\overline{\ID}$ by using the following lemma.

\begin{lem}\label{angular derivative}
Let $f\in {\mathcal U}(\lambda)$ be given by \eqref{PW-eq5}, with a
function $\phi\in B_0$. Suppose that there exists a point
$\zeta\in\T$ such that the radial limit $\phi(\zeta) =
e^{i\theta_0}$.
\begin{enumerate}
\item[{\rm (i)}] If $\cos\theta_0\in[-1,\frac{2\lambda}{1+\lambda}]$, then
$\frac{\zeta\phi'(\zeta)}{\phi(\zeta)}=1$.

\item[{\rm (ii)}]  If $\cos\theta_0\in(\frac{2\lambda}{1+\lambda},1]$, then
$$1\leq\frac{\zeta\phi'(\zeta)}{\phi(\zeta)}\leq
\frac{(1+\lambda)(1+\lambda-2\lambda\cos\theta_0)}{5\lambda^2+2\lambda+1-4\lambda(1+\lambda)\cos\theta_0}.$$
\end{enumerate}
\end{lem}
\bpf The assumption together with Julia-Wolff-Lemma (cf. \cite{Shap-93}) gives that the
angular derivative $\phi'(\zeta)$ exists and satisfies that
$$1\leq|\phi'(\zeta)|\leq+\infty.$$ On the other hand, Theorem~C 
implies that
$$0<\frac{\zeta\phi'(\zeta)}{\phi(\zeta)}\leq+\infty.$$ These facts
show that
$$1\leq\frac{\zeta\phi'(\zeta)}{\phi(\zeta)}\leq+\infty.$$

Set
$$t=\frac{\zeta\phi'(\zeta)}{\phi(\zeta)}-1\ ~\mbox{ and }~  L(\phi)(z)=\left|\left(\frac{z}{f(z)}\right)^2f'(z)-1\right|.$$
Then \eqref{PW-eq5} implies that
\[L(\phi)(z)=\left|-(1+\lambda)\left(\phi(z)-z\phi'(z)\right)+\lambda\phi(z)\left(\phi(z)-2z\phi'(z)\right)\right|<\lambda,\;\ z\in\D.\]
Let $z$ approach to $\zeta$ along the segment connecting $0$ and
$\zeta$. We obtain that
$$L(\phi)(z)\rightarrow \left|(1+\lambda)t-\lambda e^{i\theta_0}(2t+1)\right|=:R(e^{i\theta_0})(t).$$
A calculation shows that
$$
R^2(e^{i\theta_0})(t)=(1+\lambda)^2t^2+\lambda^2
(2t+1)^2-2t(2t+1)\lambda(1+\lambda)\cos\theta_0.$$

(i)\quad If $\cos\theta_0\in[-1,\frac{2\lambda}{1+\lambda}]$, then
$R^2(e^{i\theta_0})(t)$ is monotonically increasing for $t\geq0$ and
attains the minimum $\lambda^2$ at $t=0$.

If $t>0$, then $ R(e^{i\theta_0})>\lambda$, which contradicts the
fact that $L(\phi)(z)<\lambda$ for all $z\in\D$. Therefore, we have
$t=0$ and $$\frac{\zeta\phi'(\zeta)}{\phi(\zeta)}=1.$$

(ii) If $\cos\theta_0\in(\frac{2\lambda}{1+\lambda},1]$, then
$R^2(e^{i\theta_0})(t)$ is monotonically decreasing for
  $t\in[0,t_0]$ and monotonically increasing for  $t\in[t_0,\infty)$,
which implies that $R^2(e^{i\theta_0})(t)\leq\lambda^2$ when
$t\in[0,2t_0]$ and $R^2(e^{i\theta_0})(t)>\lambda^2$ when
$t\in(2t_0,\infty)$, where
\[t_0=\frac{\lambda\left((1+\lambda)\cos\theta_0-2\lambda\right)}{5\lambda^2+2\lambda+1-4\lambda(1+\lambda)\cos\theta_0}.\]
The above facts and the assumption imply that
$$1\leq\frac{\zeta\phi'(\zeta)}{\phi(\zeta)}\leq
\frac{(1+\lambda)(1+\lambda-2\lambda\cos\theta_0)}{5\lambda^2+2\lambda+1-4\lambda(1+\lambda)\cos\theta_0}.$$
This completes the proof of the lemma.
 \epf

Naturally, it follows that

\begin{thm}\label{impossible condition2}
Let $f\in {\mathcal U}(\lambda)$ be given by \eqref{PW-eq5} with a
function $\phi\in B_0$. If there exists a point $\zeta\in\T$ such
that the radial limit $\phi(\zeta)\in\T$ and
$\cos\phi(\zeta)\in[-1,\frac{2\lambda}{1+\lambda}]$. Then $\phi$
must be of the form $\phi(z)= e^{i\theta} z.$
\end{thm}

\bpf By assumption, we may write $\phi(z) = z\psi(z)$, where
$\psi(\ID)\subset\ID$ or $|\psi(z)|\equiv1$ on $\ID$. Since
$\phi(\zeta)\in\T$, we have that $\psi(\zeta)\in\T$. If
$\psi(\ID)\subset\ID$, then Theorem~C 
implies that the angular derivative $\psi'(\zeta)$ exists and
$$0<\frac{\zeta\psi'(\zeta)}{\psi(\zeta)}\leq+\infty,$$ which yields
that $$\frac{\zeta\phi'(\zeta)}{\phi(\zeta)} =1 +
\frac{\zeta\psi'(\zeta)}{\psi(\zeta)}>1.$$ This contradicts Lemma
\ref{angular derivative}. So $|\psi(z)|\equiv1$ on $\ID$. The proof
is completed.
 \epf

Since each finite Blaschke product maps $\T$ onto $\T$, by Theorem
\ref{impossible condition2}, we have the following conclusion.

\begin{cor}\label{finite Blaschke}
Let $f\in {\mathcal U}(\lambda)$ be given by \eqref{PW-eq5}, with
$\phi$ being a finite Blaschke product and $\phi(0)=0$. Then $\phi$
must be of the form $\phi(z)= e^{i\theta} z.$
\end{cor}

For $\lambda=1$, then $\frac{2\lambda}{1+\lambda}=1$. Theorem
\ref{impossible condition2} yields that

\begin{cor}\label{angular derivative of 1 case}
Let $f\in {\mathcal U}$ be given by \eqref{PW-eq5} with a function
$\phi\in B_0$. Suppose that there exists a point $\zeta\in\T$ such
that the radial limit $\phi(\zeta)\in \T$. Then $\phi$ must be of
the form $\phi(z)= e^{i\theta} z.$
\end{cor}

Since each inner function has radial limit of modulus one a.e., we
obtain that

\begin{cor}\label{inner fucntion, case 1}
Let $f\in {\mathcal U}$ be given by \eqref{PW-eq5}, with $\phi$
being an inner function and $\phi(0)=0$. Then $\phi$ must be of the
form $\phi(z)= e^{i\theta} z.$
\end{cor}

%

\begin{example}\label{ not in U} Let $r_n=1-\frac{1}{2^{4n}},$ $\theta_n=\frac{1}{2^n}$ and
$a_n=r_ne^{i\theta_n}$, where $n\geq1$. Then
$$B_1(z)=z\prod^{\infty}_{n=1}\frac{|a_n|}{a_n}\frac{a_n-z}{1-\overline{a_n}z}$$
is a Blaschke product, where
\begin{enumerate}
\item[{\rm (a)}] $B_1(z)$ and all its subproducts have unimodular radial limit at every point of $\T$;

\item[{\rm (b)}] $B_1(z)$ has finite angular derivative at every point of
$\T$ and $$1<\frac{\zeta B'_1(\zeta)}{B_1(\zeta)}<\infty,\;\
\forall\zeta\in\T;$$

\item[{\rm (c)}] $B_1$ has finite phase derivative at every point of
$\T$;

\item[{\rm (d)}]  \[\frac{z}{\left(1-B_1(z)\right)^2}\notin \mathcal{U}.\]
\end{enumerate}
\end{example}

\bpf It is obvious that
$$\sum^{\infty}_{n=1}(1-r_n)=\sum^{\infty}_{n=1}\frac{1}{2^{4n}}<\infty,$$
and thus $B_1(z)$ is a Blaschke product.

(a)\quad  We can see that $a_n\rightarrow1$ as $n\rightarrow\infty$.
So for each $e^{i\theta}\in \T\setminus\{1\}$, $\{a_n\}$ is bounded away
from such $e^{i\theta}$ and
$$\sum^{\infty}_{n=1}\frac{1-|a_n|}{|e^{i\theta}-a_n|}=\sum^{\infty}_{n=1}\frac{1}{2^{4n}}\frac{1}{|e^{i\theta}-a_n|}<\infty.$$
Since $|1-a_n|^2=(1-r_n)^2+4r_n\sin^2{\frac{\theta_n}{2}}$ and
$\frac{2}{\pi}|x|\leq\sin|x|\leq|x|$ when $|x|\leq\frac{\pi}{2},$ we
have
 \beq\label{LiPoWi*-Eqn 11 in example 4}
|1-a_n|\geq2\sqrt{r_n}\bigg|\sin\frac{\theta_n}{2}\bigg|\geq\frac{2\sqrt{r_n}}{\pi}\theta_n=\frac{2\sqrt{1-\frac{1}{2^{4n}}}}{\pi}\theta_n\geq \frac{\sqrt{15}}{2\pi}\theta_n=\frac{\sqrt{15}}{2\pi}\frac{1}{2^n}
\eeq
which yields that
$$\sum^{\infty}_{n=1}\frac{1-|a_n|}{|1-a_n|}\leq\frac{2\pi}{\sqrt{15}}\sum^{\infty}_{n=1}\frac{1}{2^{3n}}<\infty.$$
The above estimates together with the result in \cite{F} imply that
(a) is proved.
\medskip

(b)\quad We consider
$$G(B_1, e^{i\theta}):=\sum^{\infty}_{n=1}\frac{1-|a_n|}{|e^{i\theta}-a_n|^2}.$$
Since $\{a_n\}$ is bounded away from $\T\setminus\{1\}$, for each
$e^{i\theta}\in \T\setminus\{1\}$, we have
$$G(B_1, e^{i\theta})=\sum^{\infty}_{n=1}\frac{1}{2^{4n}}\frac{1}{|e^{i\theta}-a_n|^2}<\infty.$$
For $\theta=0,$ \eqref{LiPoWi*-Eqn 11 in example 4} gives
$$G(B_1, 1)\leq\frac{\pi^2}{15}\sum^{\infty}_{n=1}\frac{1}{2^{2n-2}}<\infty.$$
The above estimates together with the result in \cite{F} imply that
for each $\zeta\in\T$,
$$1<\frac{\zeta
B'_1(\zeta)}{B_1(\zeta)}=1+\sum^{\infty}_{n=1}\frac{1-|a_n|^2}{|e^{i\theta}-a_n|^2}<\infty.$$
Thus the conclusion (b) is proved.
\medskip
The conclusion (c) follows from (b). Finally, (d) follows from (b) and Lemma \ref{angular
derivative}.
 \epf

\begin{example}\label{ not in U()} Let $r_n=1-\frac{1}{2^{2n}}$, $\theta_n=\frac{1}{2^n}$,
$a_n=r_ne^{i\theta_n}$, where $n\geq1$. Then
$$B_2(z)=z\prod^{\infty}_{n=1}\frac{|a_n|}{a_n}\frac{a_n-z}{1-\overline{a_n}z}$$
is a Blaschke product, where
\begin{enumerate}
 \item[{\rm (a)}] $B_2(z)$ and all its subproducts have unimodular radial limit at every point of $\T$;

\item[{\rm (b)}] $B_2(z)$ has finite angular derivative at every point of
$\T\setminus\{1\}$ and $$\frac{B'_2(1)}{B_2(1)}=\infty;$$

\item[{\rm (c)}] $$\frac{z}{\left(1-B_2(z)\right)\left(1-\lambda
B_2(z)\right)}\notin{\mathcal U}(\lambda).$$
\end{enumerate}
\end{example}

\bpf By similar reasoning as in the proof of Example \ref{ not in
U}, we obtain that $B_2(z)$ is a Blaschke product, $B_2(z)$ and all
its subproducts have unimodular radial limit at every point of $\T$,
and $B_2(z)$ has finite angular derivative at every point of
$\T\setminus\{1\}$.

Since
$$|1-a_n|^2=(1-r_n)^2+4r_n\sin^2{\frac{\theta_n}{2}}\leq(1-r_n)^2+r_n\theta^2_n=\frac{1}{2^{2n}},
$$
it follows easily that
$$G(B_2,1)=\sum^{\infty}_{n=1}\frac{1}{2^{2n}}\frac{1}{|1-a_n|^2}\geq
\sum^{\infty}_{n=1}\frac{1}{2^{2n}}2^{2n}=\infty.$$
Finally, (c) follows from (b) and Lemma \ref{angular derivative}.
 \epf

\section{Harmonic mappings generated by functions in ${\mathcal U}_2(\lambda)$ }\label{LSW18-sec4}

A complex-valued function $F$ is called harmonic in $\ID$ if ${\rm Re}\,F$ and ${\rm Im}\,F$ are real harmonic in $\ID.$
Such function has the canonical representation $F=H+\overline{G},$ where $H$ and $G$ are analytic in $\ID$ with $G(0)=0.$
The Jacobian $J_{F}$ of $F$ is given by $J_{F}=|H'|^2-|G'|^2.$ Consequently, $F$ is locally univalent and sense-preserving if
and only if  $J_{F}(z)>0$ in $\ID$; or equivalently, $H'\neq 0$ in $\ID$ and the dilatation $\omega _F=:\omega =G'/H'$ belongs to $B$. See \cite{CS,Lewy}.

\begin{Lem}\label{close-to-convex condition}{\rm \cite{PK-15}}
Let $F=H+\overline{G}$ be a harmonic mapping normalized by
$H(0)=G(0)=H'(0)-1=0$ and $h\in{\mathcal S}$ be convex in $\ID$. If
$F$ satisfies that
$${\rm Re\,}\left\{e^{i\gamma}\frac{H'(z)}{h'(z)}\right\}>\left|\frac{G'(z)}{h'(z)}\right| ~\mbox{ for all $z\in\ID$}
$$
and for some $\gamma\in\IR$, then  $F$ is sense-preserving and close-to-convex (univalent) in $\ID$.
\end{Lem}

\begin{thm}\label{close-to-convex}
Suppose that $H\in{\mathcal U}_2(\lambda)$ and $0<\lambda\leq\sqrt{2}-1$. If $F=H+\overline{G}$ with the dilatation
$\omega(z)=\frac{G'(z)}{H'(z)}$, where $\omega(0)=0$ and
$$|\omega(z)|\leq\frac{1-2\lambda-\lambda^2}{(1+\lambda )^2}
~\mbox{ for all $z\in\ID$,}
$$
then $F$ is sense-preserving and close-to-convex (univalent) in $\ID$.
\end{thm}

\bpf By using Corollary 3.6 in \cite{OPSV}, we know that $H\in{\mathcal S}^{\star}$ for
$0<\lambda\leq\sqrt{2}-1<\frac{\sqrt{2}}{2}$. Set  $h'(z)=\frac{H(z)}{z}$. Then $h\in{\mathcal S}$ and is convex in
$\ID$. By Lemma~D, 
it suffices to prove that for all $z\in\ID$,
\be\label{comparison of real part and modulus}
{\rm Re}\left\{\frac{zH'(z)}{H(z)}\right\}>\left|\omega(z)\frac{zH'(z)}{H(z)}\right| ~\mbox{ for all $z\in\ID$}.
\ee
Since $H\in{\mathcal U}_2(\lambda)$, the representation \eqref{OPW7-eq1} for $H$ shows that
$$H(z)=\frac{z}{1+\lambda z\int_{0}^{z}\omega (t)\, dt},$$
where $\omega\in B$. Then
$$\frac{zH'(z)}{H(z)}=\frac{1-\lambda z^2\omega (z)}{1+\lambda z\int_{0}^{z}\omega (t)\, dt},$$
and \eqref{comparison of real part
and modulus} is equivalent to
$$|\omega(z)|<\frac{{\rm Re}\left\{\left(1-\lambda \overline{z^2\omega (z)}\right)\left(1+\lambda z\int_{0}^{z}\omega (t)\, dt\right)\right\}}
{\left|1-\lambda z^2\omega (z)\right| \left|1+\lambda
z\int_{0}^{z}\omega (t)\, dt\right|}=:M.
$$
Set $|z|=r<1$. Then we have
$$M\geq A(r, \lambda):=\frac{1-2\lambda r^2-\lambda^2 r^4}{(1+\lambda r^2)^2},
$$
and since
$$\min _{r\in [0,1]}  A(r, \lambda)=\frac{1-2\lambda-\lambda^2}{(1+\lambda )^2}
$$
we see that the conclusion holds since $1-2\lambda r^2-\lambda^2 r^4>0$ when
$0<\lambda\leq\sqrt{2}-1$.
 \epf

\begin{thm}\label{positive derivative}
Suppose that $H\in{\mathcal U}_2(\lambda)$ and
$0<\lambda\leq\frac{1}{2}$. If $F=H+\overline{G}$ with the
dilatation $\omega(z)=\frac{G'(z)}{H'(z)}$, where $\omega(0)=0$ and
$$|\omega(z)|\leq \sqrt{(1-\lambda^2)(1-4\lambda^2)^2}
~\mbox{ for all $|z|=r<1$,}
$$
then $F$ is sense-preserving and close-to-convex (univalent) in $\ID$.
\end{thm}

\bpf Using the method of the proof of Theorem \ref{close-to-convex},
we see that
$$H'(z)=\frac{1-\lambda z^2\omega (z)}{\left(1+\lambda z\int_{0}^{z}\omega (t)\, dt\right)^2},$$
where $\omega\in B$.  Set $|z|=r$ and
$F_\varepsilon(z)=H(z)+\varepsilon G(z)$, where $|\varepsilon|=1$.
Then
\[F'_\varepsilon(z)=\left(1+\varepsilon \omega(z)\right)H'(z)~\mbox{ and }~ |\arg F'_\varepsilon(z)|\leq \arcsin(|\omega(z)|)+3\arcsin(\lambda r^2).
\] Now it comes to prove
$$\arcsin(|\omega(z)|)+3\arcsin(\lambda r^2)\leq\frac{\pi}{2},$$
which is equivalent to  $|\omega(z)|\leq B(r, \lambda) $,
where
$$B(r, \lambda)= \cos(3\arcsin(\lambda
r^2))=\sqrt{1-\sin^2(3\arcsin(\lambda r^2))}=\sqrt{1-(3\lambda
r^2-4\lambda^3r^6)^2},$$ which is simplified to
$$B(r, \lambda)=\sqrt{(1-\lambda^2r^4)(1-4\lambda^2r^4)^2}.$$
It is a simple exercise to see that
$$\min _{r\in [0,1]}  B(r, \lambda)=\sqrt{(1-\lambda^2)(1-4\lambda^2)^2}.
$$
This completes the proof.
 \epf

\subsection*{Acknowledgements}
The work of the first author is supported by NSF of Hunan (No.
2020JJ6038), the Scientific Research Fund of Hunan Provincial
Education Department (20A070), the Science and Technology Plan
Project of Hunan Province (No. 2016TP1020), and the
Application-Oriented Characterized Disciplines, Double First-Class
University Project of Hunan Province (Xiangjiaotong [2018]469).

\end{document}